\documentclass{svmult}

\usepackage{makeidx}         
\usepackage{graphicx}        
\usepackage{multicol}        
\usepackage[bottom]{footmisc}
\usepackage{amsmath}
\usepackage{amsfonts}
\usepackage{amssymb}
\usepackage{bbm}
\usepackage{mathrsfs}
\usepackage{latexsym}

\newcommand{\levy}{\mathscr{B}}
\newcommand{\sheet}{\mathbbm{B}}
\newcommand{\sifbm}{\mathbf{B}}
\newcommand{\mpfbm}{\mathbf{B}}
\newcommand{\alphar}{\texttt{\large $\boldsymbol{\alpha}$}}
\newcommand{\fin}{$\Box$\\}

\makeindex             

\begin{document}

\title*{The Multiparameter Fractional Brownian Motion}
\author{Erick Herbin\inst{1}\and
Ely Merzbach\inst{2}}
\institute{Dassault Aviation, 78 quai Marcel Dassault, 92552 Saint-Cloud Cedex,
France \texttt{erick.herbin@dassault-aviation.fr} \and Dept. of Mathematics,
Bar Ilan University, 52900 Ramat-Gan, Israel \texttt{merzbach@macs.biu.ac.il}}
%
%



\maketitle

{AMS classification: 62\,G\,05, 60\,G\,15, 60\,G\,17, 60\,G\,18}

\begin{abstract}
We define and study  the multiparameter fractional Brownian motion. This
process is a generalization of both the classical fractional Brownian motion
and  the multiparameter Brownian motion, when the condition of independence is
relaxed. Relations with the L\'evy fractional Brownian motion and with the
fractional Brownian sheet are discussed. Different notions of
 stationarity of the increments for a multiparameter process are studied
and
 applied to the fractional property. Using  self-similarity we present a
 characterization for such processes. Finally, behavior of  the
multiparameter fractional Brownian motion along increasing paths is analysed.
\end{abstract}
 
\section{Introduction}
The aim of this paper is to give a satisfactory definition of the concept of
Multiparameter Fractional Brownian Motion (MpfBm). The definition given here is
a particular case of the Set-indexed Fractional Brownian Motion studied in
$\cite{ehem}$, but in the multiparameter case, the various stationarity properties can be compared.

In the last decade, two other definitions for the MpfBm appeared in the
literature (see $\cite{herbin}$ for a review of their properties). Both are
problematic as extensions of the classical fractional Brownian motion. In this
work, we hope to persuade the reader that our definition is natural, is the
``right" generalization of the fractional Brownian motion (fBm) and can be
applied directly to real applied problems.

\section{Definition of the MpfBm}\label{mpfbm}

In \cite{ehem}, a set-indexed extension of fractional Brownian motion was
defined and some extensions of fractal properties were established. Let
$\mathcal A$ be an indexing collection
 of compact subsets of a metric measure
space (metric $d$ and measure $m$) satisfying certain assumptions, the {\em
Set-indexed fractional Brownian motion (SifBm)} was defined as the centered
Gaussian process $\sifbm^H=\left\{\sifbm^H_U;\;U\in\mathcal{A}\right\}$ such
that
\begin{equation}\label{defsifbm}
\forall \ U,V\in\mathcal{A};\  E\left[\sifbm^H_U \sifbm^H_V\right]=\frac{1}{2}
\left[m(U)^{2H}+m(V)^{2H}-m(U\bigtriangleup V)^{2H}\right],
\end{equation}
where $0<H\leq\frac{1}{2}$ and $m$ is a measure defined on the $\sigma$-algebra
generated by~$\mathcal{A}$.

As the collection $\mathcal
A=\left\{[0,t];\;t\in\mathbf{R}^N_+\right\}\cup\left\{\emptyset\right\}$ is a
particular indexing collection, definition (\ref{defsifbm}) provides a
multiparameter process which can be seen as a multiparameter extension of
fractional Brownian motion. We get the following definition.

\begin{definition}
The {\em Multiparameter Fractional Brownian Motion (MpfBm)} is defined as the
centered Gaussian process
$\mpfbm^H=\left\{\mpfbm^H_t;\;t\in\mathbf{R}^N_+\right\}$ such that
\begin{equation}\label{defmpfbm}
\forall\ s,t\in\mathbf{R}^N_+;\  E\left[\mpfbm^H_s
\mpfbm^H_t\right]=\frac{1}{2}
\left[m([0,s])^{2H}+m([0,t])^{2H}-m([0,s]\bigtriangleup [0,t])^{2H}\right]
\end{equation}
where $m$ is a measure on $\mathbf{R}^N$ and $H\in (0,1/2]$ is called the {\em
index of similarity}.
\end{definition}

This definition looks very natural since it relies on a set-indexed process and
thus, structure of the space $\mathbf{R}^N$ is only present in indices and not
in the shape of the covariance function.

Notice that the definition of the MpfBm depends on the measure $m$.

In the particular case of $\mathbf{R}^2_+$ with the Lebesgue measure $m$, we
can explicitly give the covariance between $s=(s_1,s_2)$ and $t=(t_1,t_2)$
\begin{align*}
E\left[\mpfbm^H_s \mpfbm^H_t\right]=\frac{1}{2}
&\left[ (s_1 s_2)^{2H} + (t_1 t_2)^{2H} \right. \\
&- \left. (s_1 s_2 + t_1 t_2-2(s_1\wedge t_1)(s_2\wedge t_2))^{2H} \right].
\end{align*}

Let us notice that parameter $H$ is restricted to be in $(0,1/2]$, on the
contrary to standard fractional Brownian motion, in which $H$ is in $(0,1)$.

\begin{remark}
If the measure $m$ is absolutely continuous with respect to the Lebesgue
measure, the process $\mpfbm^H$ is almost surely null on the axis.
\end{remark}

Self-similarity is the first property of MpfBm. As a particular case of the
set-indexed fractional Brownian motion, the multiparameter process inherits
its properties. It is self-similar of index $N.H$: for all $a\in\mathbf{R}_+$,
\begin{equation*}
\left\{\mpfbm^H_{at};\;t\in\mathbf{R}^N_+\right\} \stackrel{(d)}{=}
\left\{a^{NH} \mpfbm^H_{t};\;t\in\mathbf{R}^N_+\right\}.
\end{equation*}
where $\stackrel{(d)}{=}$ denotes equality of finite dimensional distributions.

\section{Comparisons with other multiparameter extensions of fBm}\label{other}

The following two  multiparameter extensions of fractional Brownian motions are
classical. Their definitions rely on a generalization of covariance structure
of fBm based on euclidian structure of $\mathbf{R}^N$. The first definition
uses the euclidian norm and the second one uses the canonical basis of
$\mathbf{R}^N$.

\subsection{The L\'evy fractional Brownian motion}
The L\'evy fractional Brownian motion (L\'evy fBm) is defined as the mean-zero
Gaussian process $\levy^H=\left\{\levy^H_t;\;t\in\mathbf{R}^N_{+}\right\}$ such
that
\begin{equation*}
\forall s,t\in\mathbf{R}^N_{+};\quad E\left[{\levy^H_s}\;
{\levy^H_t}\right]=\frac{1}{2} \left[\|s\|^{2H}+\|t\|^{2H}-\|t-s\|^{2H}\right]
\end{equation*}
where $H\in(0,1)$.

The structure of the covariance function of $\levy^H$ provides an extension of
fractional Brownian motion where the absolute value in $\mathbf{R}_+$ is
substituted with the euclidian norm of the space $\mathbf{R}^N_+$. {}From this
 point of view, the L\'evy fBm is usually called an isotropic extension of fBm.
However, with this simple generalization, the process does not seem to be
really a multiparameter process.

This simple definition allows to state directly the self-similarity property.
For all $a\in\mathbf{R}_+$,
\begin{equation*}
\left\{\levy^H_{at};\;t\in\mathbf{R}^N_+\right\} \stackrel{(d)}{=} \left\{a^H
\levy^H_{t};\;t\in\mathbf{R}^N_+\right\}.
\end{equation*}

\subsection{The fractional Brownian sheet}
The fractional Brownian sheet is defined as the mean-zero Gaussian process
$\sheet^H=\left\{\sheet^H_t;\;t\in\mathbf{R}^N_{+}\right\}$ such that
\begin{equation*}
\forall s,t\in\mathbf{R}^N_{+};\quad E\left[{\sheet^H_s}\;
{\sheet^H_t}\right]=\frac{1}{2} \prod_{i=1}^N
\left[s_i^{2H_i}+t_i^{2H_i}-|t_i-s_i|^{2H_i}\right]
\end{equation*}
where $H=(H_1,\dots,H_N)\in(0,1)^N$.

In this definition, the euclidian structure of the space $\mathbf{R}^N$ is
strongly present in the shape of the covariance function of the fractional
Brownian sheet. Particularly, this kind of tensor product of standard
fractional Brownian motions along each direction  of the canonical basis of
$\mathbf{R}^N$ seems quite artificial and lacks  generality to be really
efficient in concrete applications.

{}From the covariance structure of fractional Brownian sheet, the
self-similarity property can be easily established. For all $a\in\mathbf{R}_+$,
\begin{equation*}
\left\{\sheet^H_{at};\;t\in\mathbf{R}^N_+\right\} \stackrel{(d)}{=}
\left\{a^{\sum_j H_j} \sheet^H_{t};\;t\in\mathbf{R}^N_+\right\}
\end{equation*}

\section{Different notions of stationarity}\label{stat}

Stationarity of increments is one of the two characteristic properties of the
classical fractional Brownian motion. In the framework of multiparameter
processes, the notion of stationarity can take different forms:

\begin{itemize}
\item Stationarity against translation
\begin{equation}\label{ST}
\forall h\in\mathbf{R}^N_{+};\quad 
\left\{ X_t - X_0;\;t\in\mathbf{R}^N_{+} \right\} \stackrel{(d)}{=} 
\left\{ X_{t+h} - X_h;\;t\in\mathbf{R}^N_{+} \right\}
\end{equation}

\item  Stationarity in the strong sense
\begin{equation}\label{SSS}
\forall g\in\mathcal{G}\left(\mathbf{R}^N\right);\quad 
\left\{ X_t - X_0;\;t\in\mathbf{R}^N_{+} \right\} \stackrel{(d)}{=} 
\left\{ X_{g(t)} - X_{g(0)};\;t\in\mathbf{R}^N_{+} \right\}
\end{equation}
where $\mathcal{G}\left(\mathbf{R}^N\right)$ is the set of rigid motions on $\mathbf{R}^N$; see \cite[p.~392]{taqqu}.
\end{itemize}

For the next definitions, one needs the notion of the increment of a process
$X$ on a rectangle $D=[s,t],$\ $s=(s_1,\dots, s_N)$ and $t=(t_1,\dots,t_N)$
where $s\prec t$ $(s_i\le t_i,$\ $i=1,\dots,N)$
$$\Delta X(D)=\sum_{r\in\{0,1\}^N}(-1)^{N-\sum_ir_i}X_{[s_i+r_i(t_i-s_i)]_i}.$$
This definition can be extended to finite unions of rectangles of $\mathbf{R}^N_+$. For $C=\bigcup_{i=1}^n D_i$, where the $D_i$' are rectangles such that
\[
\forall i,j\in\left\{1,\dots,n\right\} \quad
D_i \cap D_j \ne \emptyset \Rightarrow i=j,
\]
the increment $\Delta X(C)$ is defined by
\begin{equation*}
\Delta X(C)=\sum_{i=1}^n \Delta X(D_i).
\end{equation*} 
This definition is consistent as the previous expression is independent of the representation of $C$.

\begin{itemize}
\item Increment stationarity against translation
\begin{equation}\label{IST}
\forall h\in\mathbf{R}^N_{+};\quad 
\left\{ \Delta X_{[0,t]};\;t\in\mathbf{R}^N_{+} \right\} \stackrel{(d)}{=} \left\{ \Delta X_{[h,t+h]};\;t\in\mathbf{R}^N_{+} \right\}
\end{equation}

\item  Increment stationarity in the strong sense
\begin{equation}\label{ISSS}
\forall g\in\mathcal{G}\left(\mathbf{R}^N\right);\quad 
\left\{ \Delta X_{[0,t]};\;t\in\mathbf{R}^N_{+} \right\} \stackrel{(d)}{=} \left\{ \Delta X_{[g(0),g(t)]};\;t\in\mathbf{R}^N_{+} \right\}
\end{equation}

\item  Measure stationarity (also called $\mathcal{C}_0$-increment stationarity)
\begin{equation}\label{MS}
\forall\ t,\ \forall \tau \succ \tau'\in \mathbf{R}^N_{+};\quad
m([0,\tau])-m([0,\tau'])=m([0,t]) \Rightarrow X_t - X_0 \stackrel{(d)}{=}
X_{\tau} - X_{\tau'}
\end{equation}

\item  Increment measure stationarity

For all finite unions of rectangles $C$ and $C'$,
\begin{equation}\label{IMS}
m(C)=m(C') \Rightarrow \Delta X_{C}\stackrel{(d)}{=}\Delta X_{C'}.
\end{equation}

\end{itemize}

Notice that, among these 6 properties of stationarity, the first 4 are process
properties, but the last 2 properties are pointwise properties and depend of the chosen measure $m$.

The following result summarizes the connections between these different
stationarity properties:

\begin{proposition}\label{prop4.1} The following implications hold:

\begin{align*}
&(4) \Rightarrow (3) \Rightarrow (5);\\
&(4) \Rightarrow (6) \Rightarrow (5);\\
&(8) \Rightarrow (7).
\end{align*}

\end{proposition}

From  proposition 3.6 and theorem 4.4 in \cite{ehem}, the following can be stated:

\begin{proposition}\label{prop4.2}
The MpfBm is $\mathcal{C}_0$-increment stationary, but not increment measure stationary if $H\ne\frac{1}{2}$.
\end{proposition}

Let $\mpfbm^H$ be a MpfBm. The increment covariance between two rectangles $D$ and
$D'$, \ $E[\Delta \mpfbm^H(D)\cdot \Delta \mpfbm^H(D')]$ can be computed, but the formula is quite complicated.

In the particular case of $\mathbf{R}_+^2,$ with the Lebesgue measure and
$D=D'=(s,t],$ we get:
\begin{align*}
E[\Delta
\mpfbm^H(D)]^2&=(t_1t_2-s_1t_2)^{2H}+(t_1t_2-t_1s_2)^{2H}-(t_1t_2-s_1s_2)^{2H}\\
&-(s_1t_2+t_1s_2-2s_1s_2)^{2H}+(s_1t_2-s_1s_2)^{2H}+(t_1s_2-s_1s_2)^{2H}.
\end{align*}

We summarize stationarity properties for other definitions of multiparameter fractional Brownian motion (see \cite{herbin}, \cite{taqqu} and \cite{ehem}).

\begin{proposition}
The L\'evy fractional Brownian motion $\levy^H$ ($H\in (0,1)$) satisfies (\ref{ST}), (\ref{SSS}), (\ref{IST}), (\ref{ISSS}), and the fractional Brownian sheet $\sheet^H$ ($H\in (0,1)^N$) satisfies (\ref{IST}).
Moreover, if $\sheet^H$ has constant parameter $H$ in every axis, then it satisfies (\ref{MS}).
\end{proposition}

\section{Characterization}

In Sections \ref{mpfbm} and \ref{stat}, the multiparameter fractional Brownian
motion was shown to be self-similar and $C_0$-increment stationary. As standard
fractional Brownian motion is characterized by its two fractal properties,
self-similarity and stationarity, it is natural to wonder what are the
multiparameter processes satisfying the two properties.

As a particular case of set-indexed fractional Brownian motion, the
multiparameter fractional Brownian motion satisfies a pseudo-characterization
property.

\begin{proposition}\label{characmpfbm}
Let $X=\left\{X_t;\;t\in\mathbf{R}^N_+\right\}$ be a multiparameter process
satisfying the following two  properties:
\begin{enumerate}
\item[1.] self-similarity of index $\alpha\in (0,N/2)$, \item [2.]
$\mathcal{C}_0$-increment stationarity, for Lebesgue measure $m$.
\end{enumerate}
Then, the covariance function between $s$ and $t$ such that $s\prec t$ is
\begin{equation*}
E\left[X_s.X_t\right]=K\;
\left[m([0,s])^{2\alpha/N}+m([0,t])^{2\alpha/N}-m([0,t]\setminus
[0,s])^{2\alpha/N}\right].
\end{equation*}
\end{proposition}

\begin{proof} The result simply relies on Proposition 4.1 of \cite{ehem}, where we
consider the operation of $\mathbf{R}_+$ on
$\left\{[0,t];\;t\in\mathbf{R}^N_+\right\}$ such that
\begin{equation*}
\forall a>0, \forall t\in\mathbf{R}^N_+;\quad a.[0,t] = [0,at].
\end{equation*}
In that framework, we have
\begin{equation*}
\forall a>0, \forall t\in\mathbf{R}^N_+;\quad m(a.[0,t])=a^N m([0,t])
\end{equation*}
and then, $\mu$ is the function $a\mapsto a^N$, which is surjective.
\fin
\end{proof}

A consequence of Proposition \ref{characmpfbm} is that the fractal properties
of self-similarity and $\mathcal{C}_0$-increments stationarity prescribe the covariance
between points $s$ and $t$ that are comparable for the partial order $\prec$ of
$\mathbf{R}^N$. Since there are non ordered points, we cannot get a complete
characterization of the MpfBm by the two properties of self-similarity and
stationarity.

A natural question is then, what are the self-similar processes which are
stationary in the different definitions of Section \ref{stat}? The following
result shows that for some choice of stationarity definition, we obtain
characterization of the L\'evy fBm.

\begin{proposition}\label{characlevy}
Let $H\in (0,1)$. The L\'evy fBm  is the only Gaussian process which is
self-similar of index $H$ and stationary in the strong sense (property (\ref{SSS})).
\end{proposition}

\begin{proof} (cf. \cite[p.~393]{taqqu})

It is known (Sections \ref{other} and \ref{stat}) that the L\'evy fBm  is
self-similar and has stationary increments in the strong sense.

Conversely, let $X=\left\{X_t;\;t\in\mathbf{R}^N_+\right\}$ be a Gaussian
process such that
\begin{equation*}
\forall a\in\mathbf{R}_+;\quad \left\{X_{at};\;t\in\mathbf{R}^N_+\right\}
\stackrel{(d)}{=} \left\{a^{H} X_{t};\;t\in\mathbf{R}^N_+\right\}
\end{equation*}
and
\begin{equation*}
\forall g\in\mathcal{G}\left(\mathbf{R}^N\right);\quad \left\{ X_t -
X_0;\;t\in\mathbf{R}^N_+\right\} \stackrel{(d)}{=} \left\{ X_{g(t)} -
X_{g(0)};\;t\in\mathbf{R}^N_+\right\}
\end{equation*}

First of all, considering the canonical basis $\left(\epsilon_i\right)_{1\leq
i\leq N}$ of $\mathbf{R}^N$, and the rotation $g_u$ that maps $\epsilon_1$ onto
any unit vector $u$, the stationarity property leads to
\begin{equation*}
E\left[X_u\right]=E\left[ X_{g_u(\epsilon_1)} - X_0 \right]
=E\left[X_{\epsilon_1}\right]
\end{equation*}

For any $s$ and $t$ in $\mathbf{R}^N_+$, the self-similarity property leads to
\begin{equation*}
E\left[X_t - X_s\right]= E\left[X_{t-s}\right]-E\left[X_0\right] =\|t-s\|^{H}
E\left[X_{\epsilon_1}\right].
\end{equation*}
As we also have
\begin{equation*}
E\left[X_t - X_s\right]= \left( \|t\|^{H} - \|s\|^{H} \right)
E\left[X_{\epsilon_0}\right],
\end{equation*}
we get $E\left[X_t\right]=0$ for all $t\in\mathbf{R}^N_+$.

In the same way, we prove that for any $s$ and $t$ in $\mathbf{R}^N_+$,
\begin{equation*}
E\left[(X_t - X_s)^2\right]= E\left[(X_{t-s} - X_0)^2\right] =\|t-s\|^{2H}
E\left[X_{\epsilon_1}^2\right].
\end{equation*}
The result follows. \fin
\end{proof}

In the fractional Brownian sheet  case, several supplementary assumptions are
needed to obtain a characterization of the process. Particularly, a null value
of the process on each axis must be imposed as well as a condition of
self-similarity for each parameter, when the $N-1$ other ones are fixed (see
\cite{leger}). From that point of view, the fractional Brownian sheet has no
real motivation to be considered, although it satisfies the two properties of
stationarity and self-similarity.

\section{Projection on flows and regularity}

The notion of flow is the key to reduce  the proof of many theorems. It was
extensively studied in \cite{cime} and \cite{Ivanoff}.

\begin{definition} Let $S=[a,b]\subseteq R.$ An increasing function $f: S\to
\mathbf{R}^N_+$\linebreak $(x<y\Rightarrow f(x)\prec f(y))$ is called a {\em
flow}.
\end{definition}

The following results, proved in \cite{ehem}, give a good justification of the
definition of the MpfBm.

\begin{proposition}\label{prop7.1}
Let $\mpfbm^H$ be a MpfBm and $f$ be a flow. Then the process
$(\mpfbm^H)^f=\{\mpfbm_{f(t)}^H,\ t\in [a,b]\}$ is a time changed fractional Brownian motion.
\end{proposition}

However, in general, the projection of a multiparameter process does not
inherit its different properties.

\begin{proposition}\label{prop7.2}
Let $f$ be a flow, and $X$ be a multiparameter process.
\begin{enumerate}\item[1.] If $X$ is a L\'evy fBm, then
$(X)^f$ is a classical fractional Brownian motion iff $f(t)=\alpha t$ where
$\alpha\in\mathbf{R}_+^N.$ \item[2.] If $X$ is a fractional Brownian sheet,
then $(X)^f$ is a classical fractional Brownian motion iff $f$ is a line
parallel to one axis of $\mathbf{R}_+^N.$
\end{enumerate}

\end{proposition}

We conclude this section  by giving an interpretation of $H$ parameter.

Let us recall the definition of the two classical H\"older exponents of a
stochastic process $X$ at $t_0\in\mathbf{R}_{+}$ :
\begin{itemize}
\item the pointwise H\"older exponent
\begin{equation*}
\alphar_X(t_0)=\sup\left\{ \alpha:\;\limsup_{\rho\rightarrow 0}
\sup_{s,t\in \mathcal{B}(t_0,\rho)} \frac{|X_t-X_s|}{\rho^{\alpha}} < \infty
\right\}
\end{equation*}

\item the local H\"older exponent
\begin{equation*}
\widetilde{\alphar}_X(t_0)=\sup\left\{ \alpha:\;\limsup_{\rho\rightarrow 0}
\sup_{s,t\in \mathcal{B}(t_0,\rho)} \frac{|X_t-X_s|}{|t-s|^{\alpha}} < \infty
\right\}
\end{equation*}

\end{itemize}

\begin{corollary}
Let $\mpfbm^H$ be a multiparameter fractional Brownian motion with self-similarity index $H$.
The pointwise and local H\"older exponents of the projection $(\mpfbm^H)^f$ 
along any flows $f$ at $t_0\in [0,1]$, satisfy almost surely
\begin{align*}
\alphar_{(\mpfbm^H)^f}(t_0)&=\left\{
\begin{array}{ll}
\alpha_{\theta}(t_0).H & \textrm{if }\alpha_{\theta}(t_0)<1\\
H & \textrm{otherwise}
\end{array}
\right. \\
\widetilde{\alphar}_{(\mpfbm^H)^f}(t_0)&=\left\{
\begin{array}{ll}
\tilde{\alpha}_{\theta}(t_0).H & \textrm{if }\tilde{\alpha}_{\theta}(t_0)<1\\
H & \textrm{otherwise}
\end{array}
\right.
\end{align*}
where $\theta$ is the real function such that
$\theta(t)=m\left[f(t)\right]$ ($\forall t\in [0,1]$),
and $\alpha_{\theta}(t_0)$ (resp. $\tilde{\alpha}_{\theta}(t_0)$) is the 
pointwise (resp. local) H\"older exponent of $\theta$ at $t_0$.
 
\end{corollary}

Consequently, the $H$ parameter of the MpfBm $\mpfbm^H$ represents the regularity of the projection on any regular flow.
This fact gives a way to estimate $H$ from real data, in the frame of applications.

\section*{Acknowledgement} The authors wish to thank Prof. M. Dozzi for his helpful comments
and suggestions.

\printindex
\end{document}